\documentclass[twoside,11pt,leqno]{article}
\usepackage{amsfonts}

  \def\S{{\cal
S}} \def\H{{\cal H}}

\def\C{\texttt{C}} \def\iso{\textrm{iso}}

\def\B{B({\cal H})}

\newtheorem{df}{Definition}[section]
\newtheorem{thm}[df]{Theorem} \newtheorem{pro}[df]{Proposition}
 \newtheorem{ex}[df] {Example}
\newtheorem{rema}[df] {Remark} 
\def\sfstp{{\hskip-1em}{\bf.}{\hskip1em}}

\def\subject#1{\renewcommand{\thefootnote}{}\footnote
{AMS(MOS) subject classification (2010). Primary: {#1}}}

\def\keywords#1{\renewcommand{\thefootnote}{}\footnote
{Keywords: {#1}}}

\def\enddemo{\qed \endtrivlist} \expandafter\let\csname
enddemo*\endcsname=\enddemo

\def\qedsymbol{\ifmmode\bgroup\else$\bgroup\aftergroup$\fi
\vcenter{\hrule\hbox{\vrule
height.5em\kern.5em\vrule}\hrule}\egroup}
\def\qed{\ifmmode\else\unskip\nobreak\fi\quad\qedsymbol}

\title{\bf  On power Drazin normal and Drazin quasi-normal Hilbert space operators}
\author{\normalsize B.P.~Duggal, I.H.~Kim}
\date{}

\begin{document}

\maketitle \thispagestyle{empty} \vskip-16pt

\subject{Primary47A15, 47B15, 47B20.} \keywords{Drazin invertible operator, structure of $[(n,m)DN]$ operators, $n$-th root of normal operator, commutativity theorem  }
\footnote{The work of the second named author was supported by NRF(Korea) grant NO. NRF2019R1F1A1057574}
\begin{abstract} A Drazin invertible Hilbert space operator $T\in \B$, with Drazin inverse $T_d$, is $(n,m)$-power D-normal, $T\in [(n,m) DN]$, if $[T_d^n,T^{*m}]=T^n_dT^{*m}-T^{*m}T_d^n=0$; $T$ is $(n,m)$-power D-quasinormal, $T\in [(n,m) DQN]$, if $[T_d^n,T^{*m}T]=0$. Operators $T\in [(n,m) DN]$ have a representation $T=T_1\oplus T_0$, where $T_1$ is similar to an invertible normal operator and $T_0$ is nilpotent. Using this representation, we have a keener look at the structure of $[(n,m) DN]$ and $[(n,m) DQN]$ operators. It is seen that $T\in [(n,m) DN]$ if and only if $T\in [(n,m) DQN]$, and if $[T,X]=0$ for some operators $X\in\B$ and $T\in [(1,1) DN]$, then $[T^*_d,X]=0$. Given simply polar operators $S, T\in [(1,1) DN]$ and an operator $A=\left(\begin{array}{clcr} T&C\\0&S \end{array}\right) \in B(\H\oplus\H)$,  $A\in [(1,1) DN]$ if and only if $C$ has a representation $C=0\oplus C_{22}$.

\end{abstract}


\section {\sfstp Introduction} Let $\B$ denote the algebra of operators, i.e. bounded linear transformations, on a complex infinite dimensional Hilbert space $\H$ into itself. For $S,T\in\B$, let $[S,T]=ST-TS$ denote the commutator of $S, T$. An operator $A\in\B$ is normal if $[A^*,A]=0$. The spectral mapping theorem guarantees the existence of  normal $n$th roots of a normal operator $A\in\B$; however, normal $A$ may have other non-normal $n$th roots. If $T\in\B$ is an $n$th root of a normal operator $A\in\B$, then an application of the Fuglede theorem \cite{{Hal}, {K}} to $[T^n,T]=0$ implies $[T^n,T^*]=0$. Conversely, $[T^n,T^*]=0$ implies $T^n$ is normal. Recall, \cite{DR}, that $T\in\B$ is Drazin invertible if there exists an operator $T_d\in\B$ such that
$$
[T_d,T]=0, \ T^2_dT=T_d, \ T^{p+1}T_d=T^p
$$
for some integer $p\geq 1$. The operator $T_d$ is then the Drazin inverse of $T$ and $p$ is the Drazin index of $T$. A generalization of $[T^n,T^*]=0$ is obtained upon replacing $T$ by $T_d$: $T$ is Drazin normal , $T\in [DN]$,  if $[T^n_d,T^*]=0$ \cite{DY} and $T$ is $(n,m)$-Drazin normal, for some integer $m\geq 1$,  $T\in [(n,m) DN]$, if $[T^n_d,T^{*m}]=0$ \cite{AA}.

\

It is clear that if we let the positive integer $k$ denote the {\em least common multiple} of $n$ and $m$, $k={\rm LCM}(n,m)$, then $T\in [(n,m) DN]$ implies $T^k_d$ is normal. As an $n$th root of a normal operator, $T_d$, has a well defined structure \cite{{Em}, {Gi}, {RR}}. Add to this the fact that as a Drazin invertible operator, $T$ has a direct sum decomposition of type $T=T_1\oplus T_0$, $T_1$ invertible and $T_0$ nilpotent (of some order), and $T_d$ has a decomposition $T_d=T^{-1}_1\oplus 0$, it follows that $T_1$ is similar to a normal operator \cite{St}. Using this characterisation, we study the structure of $[(n,m) DN]$ operators in this note to prove that the (so called) class $[(n,m) DQN]$ of $(n,m)$ D-quasinormal opertaors $T$, $[T^n_d,T^{*m}T]=0$, studied by \cite{{DY}, {AA}} coincides with the class of $[(n,m) DN]$ operators. It is seen that $T\in [(n,m) DN]\wedge [(n+1,m) DN]$ (resp., $T\in [(n,m) DN]\wedge [(n,m+1) DN]$) if and only if $T\in [(k,m) DN]$ (resp., $T\in [(n,k) DN]$) for all integers $k\geq 1$; an $m$-partially isometric $[(n,m) DN]$ contraction is the direct sum of a unitary with a nilpotent; $[T,X]=0$ implies $[T_d^*,X]=0$ for $T\in [DN]$ and $X\in\B$.  More generally, if $A, B\in\B$ are such that $TA=BT$ for an operator $T\in [DN]$, and if either of the hypotheses $AT=TB$ and $T_d (A-B)=(B-A) T_d$ is satisfied, then $T^*_dA=BT^*_d$ and $AT^*_d=T^*_dB$.  Given operators $S, T\in [(n,m) DN]$, we prove a sufficient conditiion for the upper triangular operator $A=\left(\begin{array}{clcr}T&C\\0&S\end{array}\right)$ to be an $[(n,m) DN]$ operator; it is seen that this condition is necessary too in the case in which $n=m=1$, and  both $S$ and $T$ have a simple pole at $0$.

\

\section {\sfstp Results.}
Throughout the following, $S, T$ shall denote operators in $\B$, $n$ and $m$ shall denote positive integers, and $I$ shall denote the identity map. The spectrum of $T$ will be denoted by $\sigma(T)$ and $\iso\sigma(T)$ shall denote the isolated points of the spectrum of $T$. Many of the properties of $[(n,m) DN]$ operators lie on the surface. For example, $T\in [(n,m) DN]$ implies $T^k\in [(n,m) DN]$ for all integers $k\geq 1$, since
$$
(T^k)_d=T^k_d, \ [T^n_d,T^{*m}]=0\Longrightarrow [T_d^{kn}, T^{*mk}]=0.
$$
If $S,T\in [(n,m) DN]$ and $[S,T]=0=[S^*,T]$, then $(TS)_d=T_dS_d=S_dT_d=(ST)_d$,
$$
[T^n_d,T^{*m}]=0=[S^n_d,S^{*m}]\Longrightarrow [(TS)^n_d,(TS)^{*m}]=0,
$$
and this (result) in turn implies (for tensor product $T\otimes S$ of $T$ and $S$) that
$$
[(T\otimes S)^n_d,(T\otimes S)^{*m}]=[(T^n_d\otimes I)(I\otimes S^n_d),(T^{*m}\otimes I)(I\otimes S^{*m})]=0.
$$
For an understanding of some of the not so apparent structural properties of operators $T\in [(n,m) DN]\vee [(n,m) DQN]$, we start by recalling that $T$ is Drazin invertible if and only if $T$ has finite ascent and finite descent \cite {{DR}, {TL}}. Equivalently, $T$ is Drazin invertible if and only if $0\in\iso\sigma(T)$ and there exists an integer $p\geq 1$, called the Drazin index of $T$, such that
$$
\H=T^p(\H)\oplus T^{-p}(0)=\H_1\oplus\H_0,\ T=T|_{T^p(\H)}\oplus T|_{T^{-p}(0)}=T_1\oplus T_0.
$$
Here, $T_1$ is (evidently) invertible and $T_0$ is $p$-nilpotent. (In the case in which $0\notin\sigma(T)$, we allow ourselves a misuse of language and let $T^{-1}$ denote the Drazin inverse of $T$). Denoting as before the Drazin inverse of $T$ by $T_d$, $T_d$ has a direct sum representation
$$
T_d=T_1^{-1}\oplus 0\in B(\H_1\oplus\H_0)
$$
\cite[Theorem 2.2.3]{DR}. Evidently,
\begin{eqnarray*}
T\in [(n,m) DN] &\Longleftrightarrow& [T^n_d,T^{*m}]=0\\ &\Longleftrightarrow& [T_1^{-n}\oplus 0,T^{*m}_1\oplus T^{*m}_0]=0\\ &\Longleftrightarrow& [T^{-n}_1,T^{*m}_1]\oplus 0=0\\ &\Longleftrightarrow& [T^n_1,T^{*m}_1]=0.
\end{eqnarray*}
Hence:

\begin{pro}
	\label{pro10} $T\in [(n,m) DN]$ if and only if $T_1\in [(n,m) DN]$.
\end{pro}

The following theorem provides further information on the structure of $[(n,m) DN]$ operators $T$.

\begin{thm} \label{thm10}For every $T\in [(n,m) DN]$, there exists a direct sum decomposition $\H=\H_1\oplus\H_0$ of $\H$ and a decomposition $T=T|_{H_1}\oplus T|_{\H_0}=T_1\oplus T_0$ of $T$ such that $T$ is similar to the direct sum of a normal operator in $B(\H_1)$ with a nilpotent operator (of the order of the Drazin index of $T$) and $T_d$ is similar to a normal operator .
\end{thm}

\begin{demo} Assuming $p$ to be the Drazin index of $T$, define the (closed) subspaces $\H_1$ and $\H_0$ and the operators $T_1$ and $T_0$ as above. Then
$$
\H=\H_1\oplus\H_0,\ T=T_1\oplus T_0\in B(\H_1\oplus\H_0)
$$
(with $T_1$ invertible and $T_0$  $p$-nilpotent). Let $s={\rm LCM}(n,m)$. Then
	
$$
[T^n_d,T^{*m}]=0\Longrightarrow [T^n_d,T^{*m}_d]=0\Longrightarrow [T^s_d,T^{*s}_d]=0,
$$
i.e., $T_d^s$ is normal. Since
$$
T^s_d \ {\rm is \ normal}\Longleftrightarrow T^{-s}_1  \ {\rm is \ normal}\Longleftrightarrow T^s_1 \  {\rm is \ normal},
$$
it follows from \cite{St} that there exists an invertible normal operator $N_1\in B(\H_1)$ and an invertible operator $S_1\in B(\H_1)$ such that $T_1=S^{-1}_1N_1S_1$. Letting $S=S_1\oplus I|_{\H_0}$ and $N=N_1^{-1}\oplus 0$, we have $T_d=S^{-1}NS$.\end{demo}

Theorem \ref{thm10} leads to the simplification of the proofs of a number of results from \cite{{DY}, {AA}}. Postponing this exercise for the time being, we start here with the following proposition which (contrary to the claim in \cite{{AA}, {DY}}) proves that the classes $[(n,m) DN]$ and $[(n,m) DQN]$ of Hilbert space operators coincide.

\begin{pro} \label{pro11} $T\in [(n,m) DQN]\Longleftrightarrow T\in [(n,m) DN]$.
\end{pro}

\begin{demo}Following the notation above,
\begin{eqnarray*}
T\in [(n,m) DQN] &\Longleftrightarrow& [T^n_d,T^{*m}T]=0\\ &\Longleftrightarrow& [T^{-n}_1\oplus 0, (T^{*m}_1\oplus T^{*m}_0)(T_1\oplus T_0)]=0\\ &\Longleftrightarrow& [T^{-n}_1,T^{*m}_1T_1]=0\\ &\Longleftrightarrow& [T^{-n}_1,T^{*m}_1]=0\\
&\Longleftrightarrow& [T^{-n}_1\oplus 0, T^{*m}_1\oplus T^{*m}_0]=0\\ &\Longleftrightarrow& [T^n_d,T^{*m}]=0\\ &\Longleftrightarrow& T\in [(n,m) DN].
\end{eqnarray*}
This completes the proof.
\end{demo}

\begin{rema}\label{rema10} {\em Defining the invertible operator $S$ as in the proof of Theorem \ref{thm10}, it is seen that the operators $T\in [(n,m) DQN]\vee  [(n,m) DN]$ are similar to the direct sum of a normal operator with a nilpotent operator. Hence, for operators  $T\in [(n,m) DQN]\vee  [(n,m) DN]$,
both $T$ and $T^*$ satisfy Bishop-- Eschmeier-- Putinar properties $(\beta)_{\epsilon}$ and $(\beta)$. (The interested reader will find all pertinent information related to these properties, and results on operators satisfying these properties,  in references \cite{{EP}, {LN}, {D}}.) In particular, such operators $T$ are decomposable (hence have the single-valued extension property). Furthermore, because of similarity to the direct sum of a normal and a nilpotent operator, points $\lambda\in\iso\sigma(T)$ for such $T$ are poles of the resolvent of the operator: simple poles if $\lambda\neq 0$ and a pole of order $p$ at $0$. In consequence, operators $T$ satisfy most, generalized and classical, Browder and Weyl type theorems. ( See \cite{A} for information on Browder and Weyl type theorems.)}
\end{rema}

By definition, $T\in [(n,m) DN]\wedge [(n+1,m) DN]$ if and only if
\begin{eqnarray*}
& &  T_d (T_d^nT^{*m})=(T^{*m}T_d^n)T_d=(T^n_dT^{*m})T_d\\ &\Longleftrightarrow& T_1^{-(n+1)}T^{*m}_1= (T^{*m}_1T_1^{-n})T^{-1}_1=(T^{-n}_1T^{*m}_1)T^{-1}_1\\ &\Longleftrightarrow& [T^{-1}_1,T^{*m}_1]=0\Longleftrightarrow T_1\in [(1,m) DN]\\
 &\Longleftrightarrow& T\in [(1,m) DN];
 \end{eqnarray*}
again,  $T\in [(n,m) DN]\wedge [(n,m+1) DN]$ if and only if
\begin{eqnarray*}
& &   (T^{*m}T_d^n)T^*=(T^n_dT^{*m})T^*=T^{*(m+1)}T^n_d\\ &\Longleftrightarrow& T^{*m}_1T_1^{-n}T^*_1=T_1^{-n}T^{*(m+1)}_1=T^{*(m+1)}_1T_1^{-n}\\ &\Longleftrightarrow& [T^{-n}_1,T^*_1]=0\Longleftrightarrow T_1\in [(n,1) DN]\\
 &\Longleftrightarrow& T\in [(n,1) DN].
 \end{eqnarray*}
Hence:

\begin{pro}\label{pro11}  $T\in [(n,m) DN]\wedge [(n+1,m) DN]$ if and only if $T\in [(k,m) DN]$ and $T\in [(n,m) DN]\wedge [(n,m+1) DN]$ if and only if $T\in [(n,k) DN]$ for all integers $k\geq 1$.
\end{pro}

Proposition \ref{pro11} generalizes \cite[Propositions 2.5 --2.9]{AA}. We remark here that the hypotheses $T$ is injective in \cite[Proposition 2.6]{AA} and $T^*$ is injective in \cite[Proposition 2.9]{AA} are redundant.

\

An operator $A\in\B$ is an $m$-partial isometry for some integer $m\geq 1$ if $A^mA^{*m}A^m=A^m$. An invertible $m$-partial isometry is unitary. Hence, for operators $T\in [(n,m) DN]$ for which $T$ is an $m$-partial isometry, $T^m_1$ is a unitary, $\sigma(T)\subseteq \partial{\mathbb D}\cup\{0\}$ and $T=T_1\oplus T_0$, where $\partial{\mathbb D}$ denote the boundary of the unit disc in $\mathbb C$, $T_1$ is similar to a unitary operator \cite{St} and $T_0$ is nilpotent. Furthermore, since $T^{*m}_1=T^{-m}_1$,
\begin{eqnarray*}
T\in [(n,m) DN] &\Longleftrightarrow& [T^{-n}_1,T^{*m}_1]=0\\ &\Longleftrightarrow& T^{-n-m}_1T^{*m}_1=T^{*2m}_1T^{-n}_1=T^{*m}_1T^{-n-m}_1\\
&\Longleftrightarrow&  T\in [(m+n,m) DN].
\end{eqnarray*}
It is evident that an $m$-partially isometric operator $T\in [(n,m) DN]$ for $m=1$ is the direct sum of a unitary operator with a nilpotent: a similar conclusion holds for a general $m$ in the case in which $T$ is a contraction.

\

Recall that every contraction $A\in\B$ has a direct sum decomposition $A=A_u\oplus A_c$ into its unitary and cnu (=completely non-unitary) parts. $A$ is a cnu $C_{.0}$ contraction if $||A^{*n}x||\longrightarrow 0$ as $n\longrightarrow \infty$ for all $x\in\H$ \cite[Page 110]{K}. The operator $A$ is $k$-paranormal for some integer $k\geq 2$ if $||Ax||^k\leq ||A^kx||||x||^{k-1}$ for all $x\in\H$. It is known, see \cite[Page 319]{D1}, that $k$-paranormal contractions have $C_{.0}$ cnu parts.

\begin{pro}\label{pro12} If $T\in [(n,m) DN]$ is an $m$-partially isometric operator, then $T\in [(1,1) DN]$ (equivalently, $T\in [DN]$) and $T$ has a representation $T=U\oplus T_0$, where $U\in B(\H_1)$ is a unitary and $T_0\in B(\H_0)$ is a nilpotent.
\end{pro}

\begin{demo}
If $T\in [(n,m) DN]$ is  $m$-partially isometric, then (see above) $T^m_1$ is unitary. This, since $T$ is a contraction implies $T_1$ is a contraction, implies
$$||T_1x||^m\leq ||x||^m= ||T_1^mx||||x||^{m-1}$$ for all $x\in\H_1$. Consequently, $T_1$ is $m$-paranormal. Since $T_1$ has a non-trivial $C_{.0}$ cnu part forces $T^m_1$ to have a non-trivial $C_{.0}$ cnu part, we must have that $T_1$ is unitary. Hence $T=U\oplus T_0$ for some unitary $U$ and nilpotent
$T_0\in\ B(\H_0)$. Finally,
$$
T^{*2}_1=T^{-1}_1T^*_1=T^*_1T^{-1}_1\Longleftrightarrow [T^{-1}_1,T_1^*]=0\Longleftrightarrow [T_d,T^*]=0,
$$
i.e.,  $T\in[(1,1) DN]$.
\end{demo}

\vskip4pt {\bf Commutativity properties.} For operators $T\in [DN]$ (equivalently, $T\in [(1,1) DN]$), $T_d$ is normal, hence if $[T_d,A]=0$ for an operator $A\in\B$, then $[T_d^*,A]=0$ (by the Fuglede theorem \cite{{Hal},{K}}). Again, if $T\in [DN]$ is injective, then it is necessarily invertible and $T_d=T^{-1}$. Hence, $T$ is normal and if $[T,A]=0$ for some operator $A\in\B$, then $[T^*,A]=0=[T_d^*,A]$. The operator $T\in [DN]$ is in general not normal, and $TA=AT$ does not always imply $T^*A=AT^*$; however, $[T,A]=0$ and $T\in[DN]$ implies $[T^*_d,A]=0$, as the following argument shows. The operator $T\in[DN]$ has a direct sum representation $T=T_1\oplus T_0$,  $T_1$ invertible normal and $T_0$ nilpotent, and the Drazin inverse $T_d$ has a direct sum representation $T_d=T_1^{-1}\oplus 0$. Letting $A$ have the corresponding matrix representation $A=[A_{ij}]_{i,j=1}^2$, it is seen that $[T,A]=0$ forces $A_{12}=A_{21}=0$, and then
$$[T,A]=0\Longrightarrow [T_d,A]=0\Longleftrightarrow [T_1,A_{11}]=0\Longleftrightarrow [T^*_1,A_{11}]=0\Longleftrightarrow [T^*_d,A]=0.$$
This conclusion does not extend to $T\in [DN]$ such that $TA=BT$ for some operators $A, B\in\B$.
\begin{ex}\label{ex00} {\em  Define operators $T, A, B\in B({\mathbb C}^4)$ by $$ T=M\oplus N, \ A=\left(\begin{array}{clcr} A_1&A_3\\0&A_2\end{array}\right), \   B=\left(\begin{array}{clcr} B_1&B_3\\0&0\end{array}\right),$$
where $M, N, A_i (1\leq i\leq 3), B_i (i=1,3)$ are the $B({\mathbb C}^2)$ operators
\begin{eqnarray*} & & M=\left(\begin{array}{clcr}0&1\\-1&0\end{array}\right), \ N=\left(\begin{array}{clcr}0&1\\0&0\end{array}\right), \ A_1=\left(\begin{array}{clcr}1&0\\1&1\end{array}\right), \ A_2=\left(\begin{array}{clcr}0&-1\\0&0\end{array}\right), \\ & & 
\ A_3=\left(\begin{array}{clcr}0&0\\0&1\end{array}\right), \ B_1=\left(\begin{array}{clcr}1&-1\\0&1\end{array}\right), \  B_3=\left(\begin{array}{clcr}1&0\\0&0\end{array}\right).\end{eqnarray*}
Then $$T_d=M_d\oplus 0, [T_d,T^*]=0 (\Longleftrightarrow T\in[DN]), \ {\rm and} \ TA=BT,$$ but $$ T^*A\neq BT^* \ {\rm and} \ T^*_dA\neq BT^*_d.$$} \end{ex}

Additional hypotheses are required for $T\in [DN]$ and $TA=BT$ to imply $T^*_dA=BT^*_d$. The following theorem considers a couple of such hypotheses.

\begin{thm}\label{pro01} Given operators $A, B, T\in\B$ such that $AT=TB$, if $T\in [DN]$ and either of the hypotheses $BT=TA$ and $(A-B)T_d=T_d(B-A)$ is satisfied, then $T^*_dA=BT^*_d$ and $T^*_dB=AT^*_d$.
\end{thm}

\begin{demo} If $T\in [DN]$, then $T_dT^*=T^*T_d$, $T\in B(\H_1\oplus \H_0)$ has a decomposition $T=T_1\oplus T_0$, $T_1$ is invertible normal, $T_0$ is nilpotent, and $T_d=T_1^{-1}\oplus 0\in B(\H_1\oplus \H_0)$. Let $A, B\in B(\H_1\oplus \H_0)$ have the matrix representations
$$
A=[A_{ij}]_{i,j=1}^2  \
{\rm and} \ B=[B_{ij}]_{i,j=1}^2.
$$
Then $AT=TB$ implies
$$
A_{11}T_1=T_1B_{11}, \ A_{12}T_2=T_1B_{12}, \ A_{21}T_1=T_2B_{21}, \ A_{22}T_2=T_2B_{22}.
$$
Since $T^p_2=0$ for some integer $p\geq 1$ and $T_1$ is invertible
\begin{eqnarray*}
& & A_{12}T_2=T_1B_{12}\Longrightarrow T^p_1B_{12}=0\Longleftrightarrow B_{12}=0 \ {\rm and}\\ & & A_{21}T_1=T_2B_{21}\Longrightarrow A_{21}T^p_1=0\Longleftrightarrow A_{21}=0. \end{eqnarray*}

\

\noindent (a) Assume to start with that $BT=TA$. Then
\begin{eqnarray*}
& & B_{21}T_1=0\Longleftrightarrow B_{21}=0\Longrightarrow B=B_{11}\oplus B_{22},\\  & &T_1A_{12}=0\Longleftrightarrow A_{12}=0\Longrightarrow A=A_{11}\oplus A_{22},
\end{eqnarray*}
 $B_{11}T_1=T_1A_{11}$ ( and $B_{22}T_2=T_2A_{22}$). Hence, since $T_1$ is normal,
\begin{eqnarray*}
& & (A_{11}+B_{11})T_1=T_1(A_{11}+B_{11})\Longleftrightarrow (A_{11}+B_{11})T^*_1=T^*_1(A_{11}+B_{11}) \ {\rm and}\\
& & (A_{11}-B_{11})T_1=-T_1(A_{11}-B_{11})\Longleftrightarrow (A_{11}-B_{11})T_1^*=-T_1^*(A_{11}-B_{11}).
\end{eqnarray*}
Consequently,
\begin{eqnarray*}
& & A_{11}T^*_1=T^*_1B_{11}\Longleftrightarrow AT^*_d=T^*_dB \ {\rm and}\\
& & B_{11}T^*_1=T^*_1A_{11}\Longleftrightarrow BT^*_d=T^*_dA.
\end{eqnarray*}

\

\noindent (b) If instead $(A-B)T_d=T_d(B-A)$, then
$$
-B_{21}T^{-1}_1=0\Longleftrightarrow B_{21}=0, \ -T^{-1}_1A_{12}=0\Longleftrightarrow A_{12}=0
$$
 and
\begin{eqnarray*}
(A_{11}-B_{11})T_1^{-1}=-T_1^{-1}(A_{11}-B_{11}) &\Longleftrightarrow&  (A_{11}-B_{11})T_1=-T_1(A_{11}-B_{11})\\
&\Longleftrightarrow&  (A_{11}-B_{11})T_1^*=-T_1^*(A_{11}-B_{11}).
\end{eqnarray*}
Since we already have ($A_{11}T_1=T_1B_{11}\Longleftrightarrow$)  $A_{11}T^*_1=T^*_1B_{11}$, once again we have $AT^*_d=T^*_dB$ and $BT^*_d=T^*_dA$.
\end{demo}

Theorem \ref{pro01} is an improved version of \cite[Theorem 4.4]{AA}: it tells us that hypothesis (4.1) and any one of the hypotheses (4.2) and (4.3) of \cite[Theorem 4.4]{AA} guarantees the validity of the theorem. 

\

If $S, T\in [(n,m) DN]$ and ${\rm{LCM}}(n,m)=k$, then $S^k_d$ and $T^k_d$ are normal, hence $T_dA=AS_d$ for an operator $A\in\B$ implies $T^k_dA-AS^k_d=0=T^{*k}_dA-AS^{*k}_d$. This, however, does not guarantee $T^*_dA-AS^*_d=0$ (contrary to the claim made in \cite[Theorem 4.3]{AA}).

\

\begin{ex}\label{ex01} {\em For operators $S,T\in B({\mathbb C}^2)$, let
$$
S=T=\left(\begin{array}{clcr} 0&1\\-1&1\end{array}\right).
$$
Then
$$ T_d=\left(\begin{array}{clcr} 1&-1\\1&0\end{array}\right), \ T_d^2T^{*3}=T^{*3}T_d^2
$$
(so that $S=T\in [(2,3) DN]$). Since $T_d^{*3}=-I$, $T^{*3}_dA=AT^{*3}_d$ for all $A\in B(\mathbb C^2)$. If, however, we let $A=T$, then
$$
[T_d,A]=0  \  {\rm \  and}  \ T^*_dA\neq AT^*_d.
$$
Observe that $T_d=T^{-1}$, hence $TA=AT$ and  $T^*A\neq AT^*$. }
\end{ex}

The following theorem considers operators $S, T\in [(n,m) DN]$ such that $S, T$ are intertwined by a quasiaffinity (i.e., an injective operator with a dense range) to prove that $S, T$ are similar to the perturbation of a normal operator by nilpotent operators.

\begin{thm}\label{thm02} If $S, T\in\B$ are such that $X$ is a quasiaffinity, $S$ and  $T$ are $[(n,m) DN]$ operators and $SX=XT$, then there exist a normal operator $N$, nilpotent operators $S_0$ and $T_0$, and invertible operators $A, B\in\B$ such that $S=A^{-1}(N\oplus S_0)A$ and $T=B^{-1}(N\oplus T_0)B$.
\end{thm}

\begin{demo} There exist positive integers $p,q$ such that
$$
S=S_1\oplus S_0\in B(S^q(\H)\oplus S^{-q}(0)), \ T=T_1\oplus T_0\in B(T^p(\H)\oplus T^{-p}(0)),
$$
 where $S_0$ is $q$ nilpotent, $T_0$ is $p$ nilpotent, $S_1=A^{-1}_1N_1A_1$ and $T_1=A^{-1}_2N_2A_2$ for some normal operators $N_1\in B(S^q(\H))$ and $N_2\in B(T^p(\H))$, and invertible operators $A_1\in B(S^q(\H))$ and $A_2\in B(T^p(\H))$. Define the invertible operators $A, B_1\in\B$ by
 $$
 A=A_1\oplus I|_{S^{-q}(0)}, \ B_1=A_2\oplus I|_{T^{-p}(0)}.
 $$
 Then
$$
A^{-1}(N_1\oplus S_0)A_1X=XB^{-1}_1(N_2\oplus T_0)B_1\Longleftrightarrow (N_1\oplus S_0)Y=Y(N_2\oplus T_0),
$$
where we have set $AXB^{-1}_1=Y$. Evidently, $Y:T^p(\H)\oplus T^{-p}(0))\longrightarrow S^q(\H)\oplus S^{-q}(0)$ is a quasiaffinity. Let $Y$ have the matrix representation $Y=[Y_{ij}]_{i,j=1}^2$. Then, since $N_1, N_2$ are invertible and $S_0, T_0$ are nilpotent, a straightforward argument shows that
$$
Y_{12}=Y_{21}=0, \ Y=Y_{11}\oplus Y_{22}, \ Y_{11} \ {\rm and} \ Y_{22} \ {\rm {are \ quasiaffinities}}.
$$
Furthermore,
$$
S_0Y_{22}=Y_{22}T_0
$$
(so that indeed $p=q$) and
$$
N_1Y_{11}=Y_{11}N_2\Longleftrightarrow N^*_1Y_{11}=Y_{11}N^*_2.
$$
But then $N_1$ and $N_2$ are unitarily equivalent normal operators, i.e., there exists a unitary $U$ and a normal operator $N$ such that $N_1=N$ and $N_2==U^*NU$. Now define the operator $B$ by $B=UA_2\oplus I|_{TS^{-p}(0)}$. Then $S=A^{-1}(N\oplus S_0)A$ and $T=B^{-1}(N\oplus T_0)B$.
\end{demo}

If $S,T\in [(n,m) DN]$, then $S\oplus T\in [(n,m) DN]$. This fails for upper triangular operator matrices (with a non-trivial entry in the $(1,2)$-place).

\begin{ex}\label{ex02} {\em Consider operators $T, C\in B(\mathbb C^2)$ and $A\in B(\mathbb C^4)$ defined by  $T=\left(\begin{array}{clcr} 0&1\\-1&1\end{array}\right)$ (as in Example \ref{ex01}) and
$$
C=\left(\begin{array}{clcr} 0&1\\0&1\end{array}\right), \ A=\left(\begin{array}{clcr} T&C\\0&T\end{array}\right).
$$
Then $T\in [(2,3) DN]$ and $A_d=\left(\begin{array}{clcr}T_d&X\\0&T_d\end{array}\right)$, where $X=\left(\begin{array}{clcr} 0&0\\-1&0\end{array}\right)$. A simple calculation shows that $A\notin [(2,3) DN]$.}
\end{ex}

If $S, T\in\B$ are $[(n,m) DN]$ operators such that $S$ has Drazin index $q$ and $T$ has Drazin index $p$, then $S=S_1\oplus S_0\in B(S^q(\H)\oplus S^{-q}(0))$ and $ T=T_1\oplus T_0\in B(T^p(\H)\oplus T^{-p}(0))$. Let $C:S^q(\H)\oplus S^{-q}(0)\longrightarrow T^p(\H)\oplus T^{-p}(0)$ have the matrix representation
$C=[C_{ij}]_{i,j=1}^2$. Then the operator
$$
A=\left(\begin{array}{clcr} T&C\\0&S\end{array}\right),
$$
is Drazin invertible with Drazin inverse
$$
A_d=\left(\begin{array}{clcr} T_d&X\\0&S_d\end{array}\right),
$$
where $X$ is the operator
\begin{eqnarray*}
X &=& \left[\sum_{j=0}^{q-1}{T_d^{j+2}CS^j}\right](I-SS_d)+(I-TT_d)\left[\sum_{j=0}^{p-1}{T^jCS_d^{j+2}}\right]-T_dCS_d\\ &=& \left(\begin{array}{clcr} -T_1^{-1}C_{11}S^{-1}_1& \sum_{j=0}^{q-1}{T_1^{-j-2}C_{12}B^j_2}\\ \sum_{j=0}^{p-1}{T^j_2C_{21}S_1^{-j-2}} & 0 \end{array}\right)
\end{eqnarray*}
\cite[2.3.12 Theorem, Page 29]{DR}. The following theorem considers the case $n=m=p=q=1$ to give a necessary and sufficient condition for $A\in [DN]$.

\begin{thm}\label{thm03} Given operators $S, T\in\B$ such that $S, T\in [DN]$, $S$ and $T$ have Drazin index $1$ and $C: S(\H)\oplus S^{-1}(0)\longrightarrow T(\H)\oplus T^{-1}(0)$ has the matrix representation $C=[C_{ij}]_{i,j=1}^2$, a necessary and sufficient condition for the operator $A\in B(\H\oplus\H)$ to be a $[DN]$ operator is that $C=0\oplus C_{22}$.
\end{thm}

\begin{demo} If $S, T$ have Drazin index $1$, then $S=S_1\oplus 0$, $T=T_1\oplus 0$, $S_1$ and $T_1$ are normal invertible and the operator $X$ (above) has the form
$$
X=\left(\begin{array}{clcr}-T^{-1}_1C_{11}S^{-1}_1&T^{-2}_1C_{12}\\C_{21}S^{-2}_1&0\end{array}\right).
$$
Given  $S, T\in [DN]$, $A\in [DN]$ if and only if
\begin{eqnarray*}
[A_d,A^*]=0 &\Longleftrightarrow& \left(\begin{array}{clcr}T_dT^*+XC^*&XS^*\\S_dC^*&S_dS^*\end{array}\right)=\left(\begin{array}{clcr}T^*T_d&T^*X\\C^*T_d&C^*X+S^*S_d\end{array}\right)\\ &\Longleftrightarrow& S_dC^*=C^*T_d, \ XS^*=T^*X, \ XC^*=0=C^*X.
\end{eqnarray*}
The equality
\begin{eqnarray*}
S_dC^*=C^*T_d &\Longleftrightarrow&S^{-1}_1C^*_{11}=C_{11}T^{-1}_1, \ S^{-1}_1C^*_{21}=0=C_{12}^*T^{-1}_1\\ &\Longleftrightarrow& S^{-1}_1C^*_{11}=C^*_{11}T^{-1}_1, \ C_{12}=C_{21}=0;
\end{eqnarray*}
\begin{eqnarray*}
XS^*=T^*X &\Longleftrightarrow& T^{-1}_1C_{11}S^{-1}_1S^*_1=T^*_1T^{-1}_1C_{11}S^{-1}_1, \ C_{21}S^{-2}_1S^*_1=0=T^*_1T^{-2}_1C_{12}\\ &\Longleftrightarrow& C_{12}=C_{21}=0, \ T^{-1}_1C_{11}S^*_1S^{-1}_1=T^{-1}_1T^*_1C_{11}S^{-1}_1\\
&\Longleftrightarrow& C_{12}=C_{21}=0, \ C_{11}S^*_1=T^*_{11}C_{11}.
\end{eqnarray*}
Considering finally the equalities $XC^*=0=C^*X$, if $C_{12}=C_{21}=0$ and $C^*_{11}T_1=S_1C_{11}^*$, then
\begin{eqnarray*}
XC^*=0=C^*X &\Longleftrightarrow& T^{-1}C_{11}S^{-1}_1C^*_{11}=0=C^*_{11}T^{-1}_1C_{11}S^{-1}_1\\ &\Longleftrightarrow& C_{11}S^{-1}_1C_{11}^*=0=C^*_{11}T^{-1}_1C_{11}\\ &\Longleftrightarrow& T^{-1}_1C_{11}C^*_{11}=0=C^*_{11}C_{11}S^{-1}_1\\ &\Longleftrightarrow& C_{11}=0.
\end{eqnarray*}
It being straightforward to verify that $XC^*=0=C^*X$ and $C_{11}=0$ implies $C_{12}=C_{21}=0$ and $C^*_{11}T_1=S_1C^*_{11}$, it follows that a necessary and sufficient condition for $A\in [DN]$ is that $C=0\oplus C_{22}\in  B( S(\H)\oplus S^{-1}(0), T(\H)\oplus T^{-1}(0))$.
\end{demo}

The proof above, in particular our consideration of the equation $XS^*=T^*X$, exploited the fact that $S_1$ and $T_1$ are normal. Since this no longer holds for $S, T\in [(n,m) DN]$ for general $n, m>1$, the necessity of the condition $C=0\oplus C_{22}$ is not clear (for the general case). The following theorem, however, shows that this condition is sufficient. Let $S$ have Drazin index $q$, $T$ have Drazin index $p$, and let  $C\in B(S^q(\H)\oplus S^{-q}(0),T^p(\H)\oplus T^{-p}(0))$ have the direct sum decomposition $C=0\oplus C_{22}$.

\begin{thm}\label{thm04} If $S, T\in [(n,m) DN]$, then $A\in [(n,m) DN]$.
\end{thm}

\begin{demo} The hypothesis $C=0\oplus C_{22}$ forces $X=0$, and then $A^n_d=T^n_d\oplus S^n_d$. Define the operator $L$ by
$$
L=0\oplus \sum_{j=0}^{m-1}{S^{*j}_2C^*_{22}T^{*(m-1-j)}_2}.
$$
Then
\begin{eqnarray*}
A^n_dA^{*m} &=& \left(\begin{array}{clcr} T^n_dT^{*m}&0\\S^n_d L&S^n_dS^{*m}\end{array}\right)\\ &=& T_d^nT^{*m}\oplus S^n_dS^{*m}=T^{*m}T^n_d\oplus S^{*m}S^n_d\\
 &=&\left(\begin{array}{clcr}T^{*m}T^n_d&0\\LT^n_d&S^{*m}S^n_d\end{array}\right)=A^{*m}A^n_d,
 \end{eqnarray*}
 i.e., $ A\in [(n,m) DN]$.
 \end{demo}

Theorem \ref{thm04} is a generalized version of \cite[Theorem 2.7]{AA} (which contrary to the claim made by the authors does not prove the necessity of the stated conditions).


\vskip10pt \noindent\normalsize\rm B.P. Duggal, 8 Redwood Grove, London W5 4SZ, England (U.K.).\\
\noindent\normalsize \tt e-mail: bpduggal@yahoo.co.uk

\vskip6pt\noindent \noindent\normalsize\rm I. H. Kim, Department of
Mathematics, Incheon National University, Incheon,  22012, Korea.\\
\noindent\normalsize \tt e-mail: ihkim@inu.ac.kr

\end{document}